                       \documentstyle[12pt,amssymb,psfig]{amsart}
                        \numberwithin{equation}{section}
                        \theoremstyle{plain}
                         \newtheorem{thm}{Theorem}[section]
                                \newtheorem{lem}[thm]{Lemma}
                                \newtheorem{pro}[thm]{Proposition}
                                \newtheorem{cor}[thm]{Corollary}
                        \theoremstyle{definition}
                                \newtheorem{rem}[thm]{Remark}

\newcommand{\psdraw}[2]
         {\begin{array}{c} \hspace{-1.3mm}
         \raisebox{-4pt}{\psfig{figure=  #1.eps,width=  #2}}
         \hspace{-1.9mm}\end{array}}

\newcommand{\hstar}{{\frak h}^*}

\newcommand{\tr}{\text{\rm tr}}

\newcommand{\ve}{\varepsilon}
\newcommand{\Z}{{\Bbb Z}}

\newcommand{\fg}{{\frak g}}

\newcommand{\hcheck}{h^\vee}

\newcommand{\id}{{\rm id}}
\newcommand{\Hom}{{\rm Hom}}

\newcommand{\Cge}{\tilde{\cal C}(\fg,\ve)}
\newcommand{\Cgez}{\tilde{\cal C}(\fg,\ve;\zeta)}

\def\1{\hbox{\rm\rlap {1}\hskip .03in{\rm I}}} 
 \def\id{{\text {id}}} \def\tr{{\text {tr}}}
  \def\Z{{\Bbb Z}}

\baselineskip=  12pt \topskip=  12pt \def\skipaline{\vskip 12pt
plus 1pt} 

  \def\Hom{{\text {Hom}}}

  \def\exp {{\text {exp}}}

\def\End {{\text {End}}}

\def\Dim {{\text {Dim}}}
  \def\dim  {{\text {\rm dim}}}

\def\id {{\text {id}}}

\begin{document}

                \title[Quantum groups and    $G$-categories]
{Quantum groups and ribbon  $G$-categories}

                \author[T. Le \& V. Turaev]{Thang  Le  and Vladimir Turaev}
\thanks{The first author is partially supported by an NSF grant.}
\address{Dept. of Mathematics, SUNY at Buffalo, Buffalo, NY 14214, USA }
\email{letu@@math.buffalo.edu}
\address{Institut de Recherche Math\'ematique Avanc\'ee, Universit\'e
Louis Pasteur - CNRS, 7 rue Ren\'e Descartes, 67084 Strasbourg
Cedex, France}
\email{turaev@@math.u-strasbg.fr}
                \maketitle
\begin{abstract}
For a group $G$, the notion of a ribbon $G$-category was
introduced in
   \cite{Tu4} with a view towards constructing $3$-dimensional
homotopy quantum field theories
(HQFT's) with target    $K(G,1)$. We discuss here how to
 derive ribbon $G$-categories from  a  simple complex Lie algebra
$\fg$ where $G$ is the center of $\fg$. Our construction is based
on  a study of representations of  the quantum group $U_q(\fg)$
at  a root  of unity $\ve$.   Under certain assumptions on $\ve$,
the resulting $G$-categories give rise to numerical invariants
of  pairs (a closed   oriented 3-manifold  $M$, an element of $
H^1(M;G)$)  and to 3-dimensional HQFT's.
\end{abstract}

\addtocounter{section}{-1}

\section{Introduction}

In order to construct $3$-dimensional homotopy quantum field
theories (HQFT's), the second author    introduced  for a group
$G$ the notions of a ribbon $G$-category and a modular (ribbon)
$G$-category. The aim of this paper is to analyze the categories
of representations of quantum groups (at roots   of unity) from
this prospective. The role of  $G$ will be played by  the center
of the corresponding   Lie algebra.

We begin with a general theory of  ribbon  $G$-categories with
abelian $G$. According to  \cite{Tu4}, any modular  (ribbon)
$G$-category $\cal C$ gives rise to a 3-dimensional HQFT  with
target   $K(G,1)$. Such an HQFT comprises several ingredients
including    a    \lq\lq   homotopy modular functor" assigning
  modules to   pairs (a surface $\Sigma$,  an element of   $
H^1(\Sigma;G)$) and a numerical invariant  of  pairs (a closed
oriented 3-manifold  $M$, a  cohomology class  $\xi \in
H^1(M;G)$). We introduce  here a  larger class   of weakly
non-degenerate premodular $G$-categories. Similarly to
\cite{Tu4}, each  such category $\cal C$   gives rise  to  a
numerical  invariant $\tau_{\cal C} (M, \xi )$ of pairs $(M,
\xi)$ as above.
  If   ${\cal C}$ satisfies an additional assumption of regularity
and $G$ is finite  then   the standard Witten-Reshetikhin-Turaev
invariant $\tau_{\cal C} (M)$ is also   defined and splits as
$$\tau_{\cal C} (M) =   \vert G\vert ^{- (b_1(M)+1)/2} \sum_{\xi
\in H^1 (M; G)} \tau_{\cal C} (M, \xi ) $$
 where
$b_1(M) $ is the first Betti number of $M$. Note that the
invariant $\tau_{\cal C} (M, \xi )$  extends to an HQFT provided
$\cal C$ is a modular $G$-category.

Let $\fg$ be a simple complex Lie algebra and $\ve$ be a  complex
root of unity. We show that under certain assumptions on the
order of
 $\ve$, the pair $(\fg, \ve)$
gives rise   to a  premodular  ribbon $G$-category  $\cal C= \cal
C(\fg,\ve)$
  where $G $ is the center of  $\fg$.  (For the center groups of
  simple Lie algebras, see Table 1 below). The definition of $\cal C$
is based on a study of the representations of the quantum group
$U_q(\fg)$,  cf. \cite{AP,Lusztig,Le2}. We specify  conditions on
$\ve$ which ensure   that $\cal C$ is regular  so that we have
numerical invariants of 1-cohomology classes and a splitting of
the standard WRT-invariant as above. Another set of conditions
ensures that   $\cal C $ is a modular $G$-category. The resulting
3-dimensional HQFT  is however not very interesting since it
splits   as a product of a standard TQFT and a homological HQFT
(cf. Remark \ref{quq}).

The   paper consists of three sections. In Sect.\ 1 we discuss
the   theory of  premodular   $G$-categories with abelian $G$. In
Sect.\ 2 we   recall the definition of  $\tau_{\cal C} (M, \xi )$
and briefly discuss the homotopy modular functor. In   Sect.\ 3
we consider the category $\cal C(\fg,\ve)$.

The first author would like to thank  A. Beliakova, C. Blanchet,
A. Bruguieres, G. Masbaum for helpful discussions.

\section {Premodular $G$--categories}

Throughout this section, $G$ denotes an abelian group and $K$ a
  field of
characteristic zero.

\subsection{Preliminaries on monoidal
categories}
We shall use the standard notions of the theory of monoidal
categories, see  \cite{Ma}. Recall that a left duality in  a
monoidal category  $\cal C$
   associates to any object $V\in {\cal C}$ an object $V^*\in {\cal C}$ and
two morphisms $ b_V:  \1\to V\otimes V^*$ and $d_V:V^*\otimes
V\to \1 $ satisfying the   identities   $$   (\text {id}_V\otimes
d_V)   \,(b_V \otimes \text {id}_V)  =  \text {id}_V,\,\,\,
(d_V\otimes \text {id}_{V^*})   (\text {id}_{V^*}\otimes b_V)
  =  \text {id}_{V^*}.  $$
Here  $\1$ denotes the unit object of $\cal C
$ and  for simplicity we omit  the
   associativity isomorphisms   and the canonical  isomorphisms
  $ \1\otimes V\approx   V\approx   V\otimes \1$. An object of $\cal C
$ isomorphic to $\1$ is said to be {\it trivial}.

A monoidal category ${\cal C}$
is
  {\it $K$-additive} if all the $\Hom's$ in ${\cal C}$ are $K$-modules and
both the
  composition and  the tensor product of morphisms are bilinear over $K$.

We say that a $K$-additive monoidal category ${\cal C}$ {\it splits
as a disjoint union of
subcategories} $\{{\cal C}_\alpha\}$ numerated by certain $\alpha $ if:

- each ${\cal C}_\alpha$
is a full subcategory of ${\cal C}$;

- each object of ${\cal C}$ belongs to ${\cal C}_\alpha$ for a unique
$\alpha $;

- if $V\in {\cal C}_\alpha$ and $W\in {\cal C}_\beta$ with
$\alpha\neq \beta$ then $\Hom_{\cal C} (V,W)=  0$.

\subsection{Ribbon $G$-categories}  \label{ribbon} A {\it  monoidal
$G$-category over $K$} is a $K$-additive monoidal category with left duality
${\cal C}$ which splits as a disjoint union of subcategories $\{{\cal
C}_\alpha\} $ numerated
by $\alpha \in G$   such that

(i) if $V\in {\cal C}_\alpha$ and $W\in {\cal C}_\beta$ then
$V\otimes W\in {\cal C}_{\alpha+\beta}$;

(ii) if $V\in {\cal C}_\alpha$ then $V^*\in {\cal C}_{-\alpha}$.

We shall write   ${\cal C}=  \amalg_\alpha {\cal C}_\alpha$ and
call the subcategories $\{{\cal C}_\alpha\} $ of ${\cal C}$ the
{\it components} of ${\cal C}$.  The category ${\cal C}_0$
corresponding to the neutral element $0\in G$ is called the {\it
neutral component} of ${\cal C}$. Conditions (i) and (ii) show
that ${\cal C}_0$ is closed under tensor multiplication and
taking the dual object. Condition (i) implies that $\1\in {\cal
C}_0$. Thus, ${\cal C}_0$ is a    monoidal category with left
duality.

The standard notions of  braidings and twists in monoidal categories
apply in this setting without any changes.
A   braiding   (resp.\ twist)  in  a  monoidal   $G$-category ${\cal C}$
   is a system of invertible
morphisms $\{c_{V,W}:V\otimes W \to    W \otimes V\}_{V,W\in {\cal
C}} $ (resp.\ $\{\theta_V:V\to V\}_{V\in \cal C}$) satisfying
   the usual conditions,  see \cite{KRT,Tu2}.
We say that a monoidal $G$-category is {\it ribbon} if it is
ribbon in the sense of \cite{Tu2}, i.e., if it has   braiding and
twist compatible with each other and with duality.

The standard theory of ribbon categories applies to any ribbon
$G$-category $\cal C$. Suppose that $L$ is a framed oriented
$m$-component link in $S^3$ whose components are ordered, and
$X_1,\dots,X_m$ are $K$-linear combinations of objects of $\cal
C$. Then there is defined the quantum Reshetikhin-Turaev
invariant $ \langle L(X_1,\dots,X_m)\rangle \in \End_{\cal
C}(\1)$. In particular, for any object $X\in \cal C$, we have a
  dimension $\dim (X)=   \langle U(X)\rangle\in
\End_{\cal C}(\1)$ where $U$ is an oriented unknot with framing
0. For any endomorphism $f:X\to X$, we have a well-defined  trace
$\tr(f)\in \End_{\cal C}(\1)$ so that $\tr(\id_X)=  \dim (X)$.

\subsection{Premodular  $G$-categories.}
Let ${\cal C}$ be a  ribbon
$G$-category.   An object $V$ of ${\cal C}$ is   {\it simple}
  if $\End_{\cal C}
(V)=  K\,\id_V$. It is clear that an object isomorphic or dual to
a simple object is itself simple. The  assumption that $K$ is a
field  and \cite[Lemma   II.4.2.3]{Tu2} imply that    any
non-zero morphism between simple objects is an isomorphism.

We say that an object $V$ of ${\cal C}$ is {\it dominated by
simple objects}  if there is a finite  set of simple objects
$\{V_{r}\}_r$    of   ${\cal C}$   (possibly with repetitions)
and   morphisms
  $\{f_r:V_{r}\to  V,g_r:V\to  V_{r}\}_r$ such that
  $\id_V=  \sum_r f_rg_r $.  Clearly, if $V\in {\cal C}_\alpha$ then
without loss of
generality we can  assume that  $V_{r}\in {\cal C}_\alpha$ for all $r$.

We say that a ribbon
$G$-category ${\cal C}$ is {\it premodular} if
it satisfies the following three axioms:

 (1.3.1) the unit object $\1 $ is simple;

(1.3.2) for each $\alpha\in G$, the set $I_\alpha$  of the
isomorphism  classes of
simple objects of ${\cal C}_\alpha$ is finite;

(1.3.3)  for each $\alpha\in G$, any object  of ${\cal C}_\alpha$ is
dominated by
simple objects of ${\cal C}_\alpha$.

According to \cite[Remark  7.8]{Tu4} the set $G_{\cal C}=
\{\alpha\in G\,\vert\, \cal C_\alpha\neq \emptyset\}$ is a
subgroup of $G$. The theory   certainly reduces to $G_{\cal C}$.
For simplicity,  we will assume from now on  that    $G= G_{\cal
C}$.

Recall that  $U$ is an  unknot with   framing $0$. Let $U^\pm$ be
an unknot with framing $\pm 1$. For each $\alpha\in G$, consider
the   formal linear combination
$$\omega_\alpha =   \sum_{i\in I_\alpha} \dim(V_i) V_i$$
where $V_i$ is a simple object in ${\cal C}_\alpha$ representing
$i \in I_\alpha$. Define
$$\Delta_\alpha =   \langle U(\omega_\alpha)\rangle =  \sum_{i\in
I_\alpha} \dim(V_i)^2,$$
$$\Delta^\pm_\alpha =   \langle U^\pm(\omega_\alpha)\rangle=   \sum_{i\in
I_\alpha} v_i^{\pm 1}\dim(V_i)^2,$$
where the twist $\theta$ acts on $V_i$ as the scalar operator
$v_i\, \id$.  We recall here the following result
\cite[Lemma 6.6.1]{Tu4}.

\begin{lem} If $V \in {\cal C}_\beta$, then
 $ V \otimes \omega_\alpha =   \dim (V)\, \omega_{\alpha +\beta}$.
\label{product}
\end{lem}

The equality here is understood as an equality in the Verlinde
algebra of  $\cal C$.  Applying to both sides the $K$-linear
homomorphism $\dim$ from the Verlinde algebra to $K$, sending the
class  of any   object  to its dimension, we obtain that $ \dim
(V)\, \Delta_\alpha=    \dim (V)\, \Delta_{\alpha +\beta}$.
Taking  as $V$ a  simple object  in ${\cal C}_{-\alpha}$ and
using the fact that $\dim (V)\neq 0$ (see   \cite[Lemma
6.5]{Tu4}) we obtain the following corollary.

\begin{cor} For every $\alpha \in
G$,
 $$\Delta_\alpha =   \Delta_0.$$
 \label{equal}
 \end{cor}

 \begin{cor} One has that  $\omega_\alpha \otimes \omega_\beta =
 \Delta_0 \, \omega_{\alpha +\beta}$.
 \label{product2}
 \end{cor}
 This follows from Lemma  \ref{product} and Corollary
 \ref{equal}.

An important property of any premodular category is the sliding
property. Here is the version for premodular G-category.

\begin{figure}[htpb]
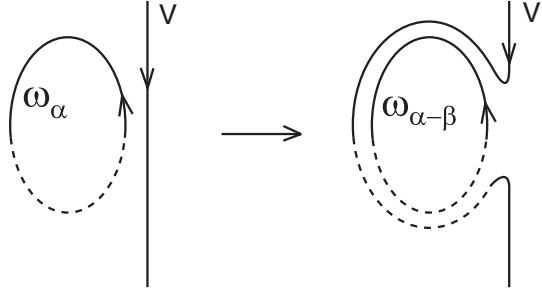

$$ \psdraw{sliding}{2.8 in} $$
\caption{Graded sliding property}\label{sliding}
\end{figure}

\begin{pro} (Graded sliding property) Let  $\cal C$  be  a   premodular $G$-category. Let $L,L'$ be framed
oriented ordered links in $S^3$ such that  $L'$ is obtained from
$L$ by sliding the second component over the first one (see Figure
\ref{sliding}). Then for every $V\in {\cal C}_\beta$, one has
\begin{equation} \langle L(\omega_\alpha, V, \dots) \rangle =   \langle
L'(\omega_{\alpha-\beta}, V, \dots) \rangle.
\label{slid}
\end{equation}
\end{pro}

The proof is the same as in the case of a premodular category,
see for example \cite{Blanchet-Beliakova}. One has to take into
account the colors of the components of tensor product
decomposition.
\begin{rem}
This proposition extends to  tangles in the obvious way. One has
to be a bit careful about colorings of a tangle. Any non-circle
component must be colored with an  object of  $\cal C$, not a
$K$-linear combination of objects    as for  circle components.
\end{rem}

\subsection{Transparent objects}
\begin{figure}[htpb]
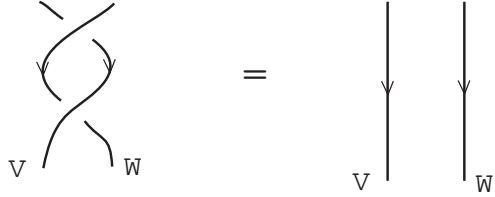

$$ \psdraw{transparent}{2.6in} $$
\caption{Transparent object}\label{transparent}
\end{figure}

Suppose that  $\cal S$ is a  set of  objects of  a ribbon
$G$-category $\cal C$.  An object $V\in \cal C$ is {\em $\cal
S$-transparent} if  for every $W\in \cal S$,
$$ c_{W,V} \, c_{V,W} =  \id.$$
This means that one can always move a string colored by $V$ past a
string colored by $W$,
 see Figure \ref{transparent}.
Let ${\cal M}_{\cal C}$ denote the set of all ${\cal
C}_0$-transparent simple objects of $\cal C$.

\begin{lem} Let $V$ be a simple object  of a
premodular $G$-category  $\cal C$  with $\Delta_0 \neq 0$. The
operator
 of the tangle in Figure \ref{S}a is non-zero for some $\alpha\in G$
if and only if $V \in {\cal M}_{\cal C}$, i.e.,    $V$ is ${\cal
C}_0$-transparent.
 \label{trans3}
 \end{lem}
\begin{figure}[htpb]
$$ \psdraw{S}{3 in} $$
\caption{}\label{S}
\end{figure}
\begin{pf} Since $V$ is simple, the operator  in Figure \ref{S} is equal to a scalar operator $t_\alpha \,
\id$ with  $t_\alpha\in K$.

Suppose $t_\alpha \neq 0$.
 Let $W\in {\cal C}_0$. Figure
\ref{transparent2}   shows that  $V \in {\cal M}_{\cal C}$; the
second equality uses the graded sliding property. (This argument
was first used by Blanchet and Beliakova  in
\cite{Blanchet-Beliakova}.)
\begin{figure}[htpb]
$$ \psdraw{transparent2}{4.2 in} $$
\caption{}\label{transparent2}
\end{figure}

Now suppose $V\in {\cal M}_{\cal C}$. It is clear that   $t_0=
\Delta_0\neq 0$. By Corollary \ref{product2}, one has
$\omega_\alpha\otimes \omega_{-\alpha}= \Delta_0 \omega_0$. The
operator of the tangle in Figure \ref{S}b is $t_\alpha\,
t_{-\alpha} \, \id$. It is also equal to (by combining the two
parallel components) $\Delta_0 \, t_0\, \id= \Delta_0^2\, \id$.
Thus $t_\alpha \, t_{-\alpha} = \Delta_0^2$, and $t_\alpha\neq 0$.
\end{pf}

\begin{rem} Under the assumption of  Lemma \ref{trans3}, if $V$ is
${\cal C}_0$-transparent, then $V$ is \lq\lq almost" transparent
with respect to any object in $\cal C$. More precisely, the same
argument as in the proof of Lemma \ref{trans3} shows that for any
$W \in {\cal C}_\alpha$, one has

\begin{figure}[htpb]
$$ \psdraw{transparent3}{2 in} $$
\caption{}\label{transparent3}
\end{figure}
\end{rem}

\subsection{Regular $G$-categories}
 A premodular
$G$-category  ${\cal C}$ is {\em weakly non-degenerate} if
 $ \Delta^+_0 \, \Delta^-_0\neq 0$. It is known (see \cite{Bru}) that
 weak non-degeneracy  implies  $\Delta_0 \neq 0$.

 A premodular $G$-category  ${\cal C}$ is {\em regular} if it is weakly non-degenerate
  and ${{\cal M}_{\cal C}}\subset {\cal C}_0$.

\begin{lem} If $\cal C$ is a regular premodular $G$-category and
$\alpha \neq 0$ then $ \Delta^{\pm}_\alpha =  0.$ \label{vanish}
\end{lem}
\begin{pf}
Consider the equality   in Figure \ref{proof4} which is obtained
by sliding the top component of the left hand side over the
bottom one.
 \begin{figure}[htpb]
$$ \psdraw{proof4}{3.6in} $$
\caption{}\label{proof4}
\end{figure}
By the regularity of ${\cal C}$, there are no ${\cal
C}_0$-transparent  simple objects  in ${\cal C}_\alpha$. Hence, by
Lemma \ref{trans3}, the right hand side on Figure \ref{proof4} is
$0$. It follows that $\Delta^+_\alpha =  0$, since
$\Delta^-_0\neq 0$. Similarly $\Delta^-_\alpha =  0$.
\end{pf}

\subsection{Modular $G$-categories} Let  $\cal C$ be a   premodular
$G$-category.
  For     $i,j \in I=  \cup_{\alpha\in G} I_\alpha$, choose
simple objects  $V_i,V_j\in {\cal C}$ representing $i,j$,
respectively, and set $$ S_{i,j}=   \langle H(V_i,V_j) \rangle=
\tr(c_{V_j,V_i}\circ   c_{V_i,V_j}:V_i\otimes V_j\to V_i\otimes
V_j)\in \End_{\cal C}(\1)=  K.$$ Here $H$ is the standard Hopf
link in $S^3$ with framing $0$.
   It follows from the properties of the trace that $S_{i,j}$ does not
depend on the
choice of $V_i, V_j$ in the isomorphism  classes $i,j$.

A premodular $G$-category $\cal C$ is {\em modular} if the
following axiom is satisfied.

({\em Non-degeneracy axiom})  The (finite) square matrix
$[S_{i,j}]_{i,j\in I_0} $ is invertible over $K$.

  It follows from this axiom that  the neutral
component ${\cal C}_0$ of ${\cal C}$   is a modular  ribbon
category in the usual, ungraded sense (see \cite{Tu2}). It is
known that the non-degeneracy   implies  the weak non-degeneracy.

A  modular ribbon $G$-category  $\cal C$   may be non-modular  in
the   ungraded sense because the set $I=  \cup_{\alpha\in G}
I_\alpha$ of  the isomorphism   classes of  simple objects in
$\cal C$ may be infinite, or because the full $S$-matrix
$[S_{i,j}]_{i,j\in I} $ may   be  non-invertible.

\begin{pro} Suppose that $G$ is finite and a modular G-category $\cal
C$ is regular. Then $\cal C$ is a modular category    in the
  ungraded sense.
\label{ttt}
\end{pro}

\begin{pf} According to  \cite {Bru}, \cite{Blanchet-Beliakova},  a
premodular category $\cal C$  is modular if and only if it has no
non-trivial $\cal C$-transparent simple objects and $\sum_{i\in
I} (\dim(V_i))^2\neq 0$. If $V$ is a $\cal C$-transparent simple
object of $\cal C$ then   $V\in {{\cal M}_{\cal C}} \subset {\cal
C}_0$. But since ${\cal C}_0$ is modular, $V$ is a trivial
object. By Corollary   \ref{equal}, $\sum_{i\in I} (\dim(V_i))^2
=  \vert G\vert \,\Delta_0
 \neq 0$, since
${\cal C}_0$ is modular.
\end{pf}

\section{Invariants of 3-manifolds and  HQFT's}

\skipaline \noindent {\bf 2.1.  Invariants of  3-dimensional
$G$-manifolds.} Fix an abelian group  $G$. Let $\cal C $ be a
weakly non-degenerate premodular ribbon $G$-category   over   a
field of zero characteristic $K$. We explain here following
\cite{Tu4} that  $\cal C $ gives rise to  a topological
invariant  of 1-dimensional cohomology classes of 3-manifolds
with coefficients in $G$.

  Fix an element ${\cal D}\in K$ such that
  ${\cal D}^2=  \Delta_0$.
Let $M$ be a closed connected oriented 3-dimensional
 manifold and $\xi \in H^1(M;G)$. Present $M$ as the result of
surgery on $S^3$ along a framed
 oriented  link $L=  L_1\cup... \cup L_m$.
Recall that $M$ is obtained by gluing $m$
  solid tori to the exterior   of $L$ in $S^3$.  This allows us to
consider for $n=  1,...,m$,   the value, $\alpha_n\in G$,  of
$\xi$ on the meridian of  $L_n$  and provide    $L_n$ with the
color $\omega_{\alpha_n}=  \sum_{i\in I_{\alpha_n}} \dim (V_i)
V_i$.
  Let $\sigma_+$ (resp. $ \sigma_-$) be  the number of positive
(resp. negative) squares in the diagonal decomposition of the
intersection form $H_2(W_L) \times H_2(W_L) \to \Bbb Z$ where
$W_L$ is the compact oriented   4-manifold  bounded by $M$ and
obtained  from the 4-ball $B^4$ by attaching  2-handles along
tubular neighborhoods of the components of $L$ in $S^3=  \partial
B^4$.
  Set
$$\tau_{\cal C,\cal D}(M,\xi) =   {\cal D}^{-b_1(M)-1}(\Delta_0^-)^{-\sigma_
-} (\Delta_0^+)^{-\sigma_+} \, \langle
L(\omega_{\alpha_1},\dots,\omega_{\alpha_m})\rangle \in K,$$ where
$b_1(M)=  m -\sigma_+-\sigma_-$ is the first Betti number of $M$.

It follows from \cite[Theorem 7.3]{Tu4} that
   $\tau_{\cal C, \cal D}(M,\xi)$ is a homeomorphism
invariant of   the pair $(M,\xi)$ (in \cite{Tu4} this is stated
for modular $G$-categories, but the proof remains true for
weakly non-degenerate premodular categories).
 The factor ${\cal D}^{-b_1(M)-1}$
appears  here for normalization purposes.

As  in the standard theory, the invariant  $\tau_{\cal C,\cal D
}(M,\xi)$ generalizes to an invariant of triples $(M,\Omega,
\xi)$ where $M$ is as above, $\Omega$ is a colored ribbon graph
in $M$ and $\xi \in H^1(M\backslash \Omega;G)$. Here a coloring is
understood as the usual coloring of $\Omega$ over $\cal C$  such
that the color of every 1-stratum $t$ of $\Omega$ is an object
in  ${\cal C}_{\alpha(t)}$ where $\alpha(t) \in G$ is the value
of $\xi$ on the meridian of $t$.  This notion of a coloring applies
in particular to framed oriented links in $M$ so that we obtain
an isotopy invariant of triples $(M,L, \xi)$ where $M$ is as
above, $L$ is a colored framed oriented link in $M$ and $\xi\in
H^1(M\backslash L;G)$.

Suppose now that  $G$ is finite and $\cal C$ is regular. Set
$\omega   =   \sum_{\alpha \in G }\omega_\alpha$. By Lemma
\ref{vanish},
$$\langle U^\pm (\omega)\rangle =   \langle U^\pm (\omega_0)\rangle
=  \Delta^\pm_0 \neq 0.$$ We can therefore consider the standard
Witten-Reshetikhin-Turaev   invariant $\tau_{\cal C, \tilde {\cal
D} }(M)$ of $M$ where $\tilde {\cal D}=   \vert G\vert   ^{1/2}
\cal D$.

\begin{pro} If $G$ is finite and  $\cal C$ is a regular premodular
$G$-category,   then
for any closed connected oriented 3-manifold $M$,
$$\tau_{\cal C,  \tilde {\cal D}}(M)=   \vert G\vert ^{- (b_1(M)+1)/2}
\sum_{\xi \in H^1 (M; G)} \tau_{{\cal C}, \cal D} (M, \xi ).$$

\end{pro}

\begin{pf}   The proof in
 \cite[section 5]{Blanchet} can be applied to get the result.
\end{pf}

\subsection {The homotopy modular functor.}  A modular $G$-category
$\cal C$ (not necessarily regular)  gives rise to a
$3$-dimensional  HQFT, see  \cite{Tu4}. In particular, the
homotopy modular functor $\cal T_{\cal C}$ associated with $\cal
C$ assigns $K$-modules to so-called extended $G$-surfaces and
assigns $K$-linear isomorphisms of these modules to weak
homeomorphisms of such surfaces. We recall here the definition of
an extended $G$-surface and the corresponding $K$-module in our
abelian case.

    Let
$\Upsilon$ be a    closed   oriented   surface. A point   $p\in
\Upsilon$  is {\it marked} if is   equipped with a sign
$\varepsilon_p=  \pm 1$ and a   tangent direction, i.e., a ray
${\Bbb R}_{+} v$ where $v$ is a non-zero tangent vector at $p$. A
{\it marking} of $\Upsilon$ is a finite (possibly void) set of
distinct marked points   $P\subset \Upsilon$.     A {\it
$G$-marking} of $\Upsilon$  is a marking $P\subset \Upsilon$
endowed with a   cohomology class $\xi \in H^1(\Upsilon\backslash
P ;G)$. A $G$-marking   $P\subset \Upsilon$ is   {\it colored} if
it is equipped  with a function $x$  which assigns to every
point   $p\in  P$   an object $x_p\in {\cal C}_{\xi (\mu_p)} $
where $\mu_p$ is a small loop in   $\Upsilon \backslash P$
  encircling
$p$ in the
  direction   induced by
the orientation of $\Upsilon$ if $ {\varepsilon_p}=  +1$ and in
the opposite direction otherwise.

An {\it extended $G$-surface} comprises a   closed oriented
surface $\Upsilon$, a colored $G$-marking   $P\subset \Upsilon$,
and a Lagrangian space $\lambda \subset H_1(\Upsilon;\Bbb R)$.
The corresponding $K$-module is  defined  as follows. Assume
that   $\Upsilon$ is a connected surface of genus $n$ and $P= \{
p_1,...,p_m\}\subset  \Upsilon$. For $r=  1,...,m$, set
$\varepsilon_r=  {\varepsilon_{p_r}}=  \pm 1, \mu_r=  \mu_{p_r}$
and $x_r=  x_{p_r} \in \cal C$. The group $\pi_1(\Upsilon
\backslash P)$  is   generated by   the homotopy classes of loops
$\mu'_{1},...,\mu'_{m}$ homotopic to $\mu_1,..., \mu_m$,
respectively,  and  by $ 2n$ elements $ a_1,b_1,..., a_n,b_n$
subject to the only relation
$$(\mu'_{1})^{\varepsilon_1}...\,(\mu'_{m})^{\varepsilon_m}\,
[a_1,b_1]\, ... \,[a_n,b_n]=  1  $$  where $[a ,b]=   ab
a^{-1}b^{-1}$. Then
$${\cal T}_{\cal C} (\Upsilon,  P, \xi, x, \lambda) $$
$$=  \bigoplus_{i_1\in I_{\xi(a_1) },
=2E..,i_n\in I_{ \xi(a_n)}}\Hom_{\cal C} \left (\1,
(x_1)^{\varepsilon_1}\otimes \dots \otimes (x_m)^{\varepsilon_m}
\otimes \bigotimes_{s=  1}^n \left   (V_{i_s} \otimes
  V_{i_s}^*\right ) \right)$$ where for an object $x \in \cal C$
we set $x^1=  x$ and $x^{-1}=  x^*$. Observe that the  dimension
of ${\cal T}(\Upsilon,  P, \xi, x, \lambda) $  over $K$ does not
depend on the values of $\xi$ on $b_1,...,b_n$. This implies that
this dimension   does not depend on the choice of $\xi$. For
$\xi=  0$, this dimension is computed by the   Verlinde formula.
In particular, if   $P=  \emptyset$ we obtain
$$\Dim \, {\cal T}_{\cal C}(\Upsilon,  P, \xi, x, \lambda)  =   \Dim \,
{\cal T}_{\cal C}(\Upsilon,  P, 0, x, \lambda)=   {\cal D}^{2n-2}
\sum_{i\in I_0}  (\dim ( V_i))^{2-2n}.$$ The same formula can be
deduced from the generalised Verlinde formula  for HQFT's given
in \cite[Sect. 4.11]{Tu3}.

\begin{rem}
Since we  restrict ourselves here to  abelian $G$, the homotopy
classes of maps from a manifold $M$ to the Eilenberg-MacLane space
$ K(G,1)$ are   classified   by elements of $H^1(M;G)$.  This
allows   one to  formulate HQFT's  with target $ K(G,1)$   in
terms of 1-dimensional cohomology classes.   We call the  HQFT's
arising from  modular $G$-categories with abelian $G$ and
trivial  crossed structure (as everywhere in this  paper)  {\it
abelian}.  The abelian HQFT's are  simpler than the general
HQFT's. One of simplifications   is that in the abelian case we
can forget about  base points   (in the general case one has
 to
work in the pointed category). Also, in the abelian case  the
canonical action of $G$ on the spaces of conformal blocks ${\cal
T}(\Upsilon,  P, \xi, x, \lambda) $ defined in \cite[Sect.
10.3]{Tu4} is    trivial.
\end{rem}

 \section{Premodular $G$-categories derived from quantum groups}

\subsection{Fusion category associated with quantum groups at
roots of unity}
We will   briefly recall the fusion categories of quantum groups
at roots of unity, mainly to fix notation. We refer the reader to
\cite{AP,GK,Kirillov} for details.

Let $\frak g$ be a finite-dimensional simple Lie algebra over
$\Bbb C$. Fix a Cartan subalgebra   and a system of basis roots
$\alpha_1,\dots,\alpha_\ell$. Choose the inner product on $\hstar
 =   {\Bbb R}$-span$\langle\alpha_1,\dots,\alpha_\ell\rangle$ such
that the square length of any {\em short} root is 2. Let $G= X/Y$
be the quotient of the weight lattice $X$ by the root lattice
$Y$. It is known that $G$ is isomorphic to the center of the
simply-connected Lie group associated with  $\frak g$, and $|G|$
is equal to the determinant of the Cartan matrix.  Let $d$ be the
maximal absolute value of  the non-diagonal entries of the Cartan
matrix. Let $h$ denote  the Coxeter number of $\fg$ and $ h^\vee$
 the dual Coxeter number of $\fg$. By $D$ we denote   the
smallest  positive integer such that $D(x|y) \in \Bbb Z$ for  all
$x,y\in X$. The data associated with simple Lie algebras is given
in Table 1.

\begin{table}[ht]
\begin{center}
\begin{tabular}{|c|c|c|c|c|c|c|c|c|c|c|c|}
\hline
     & $A_\ell$ & $B_\ell$    &  $B_\ell$  & $C_\ell$   & $D_\ell$  &
$D_\ell$  & $E_6$  &  $E_7$  &  $E_8$  & $F_4$  & $G_2$ \\
     &          & $\ell$ odd  & $\ell$ even &            & $\ell$
odd&$\ell$ even &        &         &         &        &   \\
\hline $d$  & $1$      & $2$         & 2           &  2         &
1        &   1        &   1    &    1    &    1    &   2    & 3
\\ \hline $D$  & $\ell+1$ & 2           &  1          & 1
&   4       &  2         &   3    &   2     &   1     & 1    &
1   \\ \hline $G$  &$\Z_{\ell+1}$& $\Z_2$   & $\Z_2$      &
$\Z_2$     & $\Z_4$
&$\Z_2\times\Z_2$&$\Z_3$&$\Z_2$ &   1     &  1     & 1     \\
\hline $h$  & $\ell+1$ & $2\ell$     & $2\ell$     & $2\ell$
&$2\ell-2$  &$2\ell-2$   & 12     &  18     & 30      &   12   &
6     \\ \hline $ h^\vee$
     & $\ell +1$& $2\ell-1$   & $2\ell-1$   & $\ell+1$   &$2\ell-2$
&$2\ell-2$   & 12     &  18     & 30      &   9    &   4  \\ \hline
\end{tabular}
\end{center}
\caption{}
\end{table}

The quantum group $U_v({\frak g})$, as defined in
\cite{AP,Lusztig}, is a Hopf algebra over ${\Bbb Q}(v)$, with $v$
a formal variable. There is an integral version
 of
$U_v(\fg)$, defined over the ring ${\Bbb Z}[v,v^{-1}]$,
introduced by Lusztig. For $\ve \in \Bbb C$,  let $U_\ve(\fg)$ be
the Hopf algebra over $\Bbb C$   obtained by tensoring the
integral version of $U_v({\frak g})$ with $\Bbb C$, where $\Bbb
C$ is   considered as a ${\Bbb Z}[v,v^{-1}]$-module by $v\mapsto
\ve$ (see \cite{AP,Lusztig}).

Let $\ve$ be a root of unity. The fusion category $\Cge$ is the
quotient of the category of all tilting $U_\ve(\fg)$-modules by
negligible modules and   negligible morphisms (see \cite{AP},
there $\Cge$ is denoted by ${\cal C}^-$). The category $\Cge$ is
a semisimple monoidal $\Bbb C$-abelian category with duality. The
 isomorphism classes of simple objects in $\Cge$  are parametrized by
the  dominant weights $\mu$ such that $\mu+\rho \in C_\ve$ where
$\rho$ is the half-sum of all positive roots and $C_\ve$ is
defined as follows. There are two cases depending on the  order
$r$ of  the root of unity $\ve^2$.

\label{cases}

 Case 1: The numbers $r$   and $d$ are
co-prime  and $r>h$. Then
$$ C_\ve =   \{ x \in C \, \vert \, (x|\alpha_0) < r\} $$
where $C$ is the Weyl chamber and $\alpha_0$ is the short highest
root, i.e., the only root in $C$ with square length 2.

Case 2: $r$ is divisible by $d$ and $r/d >
h^\vee$. Then
$$ C_\ve =   \{ x \in C\, \vert\, (x|\beta_0) < r\} $$
where $\beta_0$ is the long highest root, i.e., the only root in
$C$ with square length $2d$.

Note that when $d=  1$, the two cases overlap and the definitions
of $   C_\ve$ agree. In \cite{Le1}, $C_\ve$ is denoted by $C'_r$.

For  a dominant weight $\mu$ such that $\mu+\rho \in C_\ve$ there
is a simple $U_\ve(\fg)$-module $\Lambda_\mu\in \Cge$, known as
the Weyl module,   which is a deformation of the corresponding
classical $\fg$-module. It  represents  the isomorphism class
numerated by $\mu$. The decomposition of a tensor product
$\Lambda_\mu\otimes \Lambda_\nu$ is described in \cite{AP} (the
tensor product is denoted there  by \underline{$\otimes$}, but we
will use the usual notation). For our   purposes, we notice that
if $\Lambda_\gamma$ appears as a summand in $\Lambda_\mu\otimes
\Lambda_\nu$, then $\mu+\nu -w\cdot \gamma \in Y$ for some $w\in
W_\ve$ where $W_\ve$   is the group of affine transformations of
$\hstar$ generated by the reflections in   the walls of the
simplex $\bar C_\ve$, the topological closure of $C_\ve$, and
the dot action is   the one shifted by $\rho$ so that  $w\cdot x
=   w(x +\rho)- \rho$. Note that $\bar   C_\ve$ is a fundamental
domain of $W_\ve$. For further use, we state here a simple lemma.

\begin{lem} For every $x\in X$ and $w\in W_\ve$, one has  $x-w\cdot  x
\in Y$. \label{root}
\end{lem}
\begin{pf} It is known that $W_\ve =   W \ltimes Q$, where $W$ is
the usual Weyl group, and $Q$ is the $\Bbb Z$-lattice spanned by
the long roots. (One considers $Q$ as a group of translations).
Obviously, $Q\subset Y$. If $w \in Q$, the statement is trivial.
If $w\in W$, the statement is well-known in the theory of simple
Lie algebras.
\end{pf}

\subsection{Braiding and twist in $\Cge$} There are a braiding and a twist
in $\Cge$ which make  $\Cge$  a ribbon category. They depend on a
choice of a complex root of unity  $\zeta$ such that $\zeta^D=
\ve$, where $D$ is as in Table 1. Let us denote the resulting
ribbon category by $\Cgez$. The tensor structure in $\Cgez$  and
the corresponding  Verlinde   algebra do not depend on the choice
of  $\zeta$.

The ribbon category   $\Cgez$  is premodular, since it has only a
finite number of isomorphism classes of simple objects. It is
also Hermitian, see \cite{Kirillov}, so that  $\dim (V)\in  \Bbb
R$   for every $V\in \Cgez$. (This can also be  deduced  from the
explicit formulas for the quantum dimensions of the Weyl
modules). Thus for every simple object $x$, one has $(\dim x)^2$
is a positive real number (since $\dim(x)\neq 0$).

In Case 2 if $\varepsilon=  e^{\pi i/r}$, then this Hermitian
structure is known to be  unitary, see \cite{Wenzl}. Then $\dim
(V)>0$   for every $V\in \Cgez$. Again this can also be deduced
from the explicit formulas for the quantum dimensions of the Weyl
modules.

It is clear that  $\zeta^{2Dr}=  1$ where $r$ is the order of
$\ve^2$.   When $r$ is divisible by $d$, $r/d > h^\vee$,   and
$\zeta$ has  order $2Dr$, the category $\Cgez$ is modular. These
are the well-known, and so far the main, examples of modular
category.

\subsection{$\Cgez$ as a  premodular $G$-category}
Let $\pi: X \to G=   X/Y$ be the projection. For $\alpha \in G$,
let $I_\alpha$ be the set of dominant weights $\mu$ such that
$\mu+\rho\in C_\ve$ and $\pi(\mu)=  \alpha$. Let ${\cal
C}_\alpha=   {\cal C}_\alpha (\fg,\ve;\zeta)$ be the set of all
objects in $\Cgez$ which are direct sums of $\Lambda_\mu$ with
$\mu\in I_\alpha$. Set ${\cal C} =   \sqcup _{\alpha\in G}{\cal
C}_\alpha$. We consider $\cal C$ as a full subcategory of $\Cgez$.

\begin{pro} The category ${\cal C} =
\sqcup _{\alpha\in G}{\cal C}_\alpha$ is a premodular
$G$-category.
\end{pro}
\begin{pf} One only needs to check (i) and (ii) of  Sect. \ref{ribbon}.

(i) Let $\mu\in I_\alpha$,  $\nu\in I_\beta$, and
$\Lambda_\lambda$   be  a direct summand of $\Lambda_\mu\otimes
\Lambda_\nu$. Then $\mu+\nu-   w\cdot\lambda \in Y$ for some
$w\in W_\ve$.  Lemma \ref{root} shows that $\pi(\lambda) =
\pi(\mu)+\pi(\nu)$, or, in other words, $\lambda \in I_{\alpha
+\beta}$.

(ii) The dual of $\Lambda_\mu$ is $\Lambda_{-w_0(\mu)}$, where
$w_0$ is the longest element of the Weyl group. Again Lemma
\ref{root} shows that $-w_0(\mu) \in I_{-\alpha}$.
\end{pf}

If all the colors of a framed link $L$ are simple $\fg$-modules
with highest weights in the root lattice, then the quantum
invariant of $L$ is in ${\Bbb Z}[v^2, v^{-2}]$, according to the
integrality, see \cite{Le1}. Hence the fact  that $\cal C$ is a
modular, or weakly non-degenerate, $G$-category does not depend
on the choice of the $D$-th root $\zeta$ of $\ve$; it totally
depends on the order $r$ of $\ve^2$.  Note that the fact $\Cgez$
is a (non-graded) modular category does depend on the choice of
$\zeta$.

\subsection{Modular $G$-categories} The following proposition shows
that in Case 1 of Sect. \ref{cases}, the
 category  ${\cal
C}$ constructed above is a modular $G$-category, at least under
the  assumption that $r$  is    co-prime with  $|G|$.

\begin{pro} Suppose that the order $r$ of $\ve^2$ is co-prime with $d
\,|G|$ and $r > h$. Then
  $\cal C$ is a modular  $G$-category. \label{ref11}
  \end{pro}
\begin{pf} We need only to prove that the $S$-matrix of the
neutral component ${\cal
C}_0$ is invertible. This was established in \cite[Theorem 3.3]{Le2}.
\end{pf}

The following proposition shows that $\cal C$ can be $G$-modular
even when  $r$ is not co-prime with $d |G|$.

\begin{pro}\label{ref22}
Suppose $\fg$ is a Lie algebra of series $C_\ell$ with odd
$\ell$. Assume that the order $r$ of $\ve^2$ is even but not
divisible by 4 and $r
> d\hcheck$.  Then $\cal C$ is a modular $G$-category. (Note that
for $C_\ell$, one has $d=2, D=1$, and $\hcheck = \ell+1$.)
\end{pro}

A proof will be given later.

\begin{rem}   In the cases of Propositions \ref{ref11} and \ref{ref22}, it can be shown that
 the   category $\cal C $ is the product of its  neutral component and a modular
category associated with  the center group $G$ (see the
definition in \cite{Sawin}, see also \cite{Le2}). The
corresponding invariant of  a 1-cohomology class on a 3-manifold
$M$  is then the product of the invariant of the cohomology class
$0$ with an invariant depending only on $H_1(M;\Bbb Z)$ and the
linking form. The theory in this case is rather trivial.
\label{quq}
\end{rem}

\subsection{Regularity of ${\cal C}$} In   Case 2 of Sect. \ref{cases}, the
premodular $G$-category  ${\cal C}$ constructed above turns out
to be regular, at least under a few further  assumptions on
$\zeta$ and $r$.

\begin{pro} Suppose that $r$ is divisible by $d$, $r/d >  h^\vee$,
and $\zeta$ is a root of unity of order $2Dr$.  Assume that the
number $k:= r/d - \hcheck$  satisfies:

\noindent {\em (*)} \hspace{1mm} $k\ell \vdots 2|G|$ for Lie
series $A,C$; $k \vdots |G|$ for Lie series $B,E,F$; and $k\vdots
2, k\ell \vdots 2|G|$ for Lie series D,

\smallskip \noindent  where $\ell$ is the   rank of $\fg$, and
$a\vdots b$ means $a$ is divisible by $b$. Then the premodular
$G$-category $\cal C$ is  regular. \label{pre}
\end{pro}

When $\fg=   sl_n$, the condition on $r$ can be rephrased as
$r>n$  and  $(r-n)(n-1)$ is divisible by $2n$. Thus we   recover
the $sl_2$-case of \cite{KM,Tu1} and the $sl_n$-case of
\cite{Murakami,Blanchet}.

 Proposition \ref{pre} will be proven below.

 \begin{rem} One may ask when the $G$-category $\cal C$ is
 weakly-nondegenerate. The answer is only when the order $r$ of
 $\ve^2$ is as one described in  Propositions \ref{ref11}, \ref{ref22},
 \ref{pre}, or one in  the following additional cases: for the Lie
 algebra $D_\ell$ with odd $\ell$ the number $r$ must satisfies $
 r \equiv 2 \pmod 4$; for the Lie algebra $A_\ell$, one must have
 $k \ell(\ell+1) \vdots 2 s^2$, where $s = (k,\ell+1)$.
\end{rem}

\subsection{Transparent objects} Our aim here is to study
 the set ${\cal M}_{\cal C}$ of    ${\cal C}_0$-transparent
objects of  ${\cal C} $. To cover  both Cases  1 and 2  of
Sect.   \ref{cases}    we introduce a group $G'$ equal to $G$ in
Case 1 and equal to $Hom(G, \Bbb Q/\Bbb Z)$ in Case 2. More
precisely, set $G'=   X'/Y'$ where
 $X'=  X$, $Y'=  Y$ in Case 1 and
 $X'=   Y^*, Y'=  X^*$ are the $\Bbb Z$-dual lattices in Case 2.
 In both cases one has $X' \subset Y^*$.

The group $G'$ acts on $\bar C_\ve$ as follows. Let $g\in G'$
with a   lift $\tilde g\in X'$  and let $\mu\in \bar C_\ve$. The
element $r \tilde g + \mu$ may not lie in the simplex $\bar
C_\ve$ anymore, but there is  $w\in W_\ve$ which maps $r \tilde g
+ \mu$ into $\bar C_\ve$. Set $g(\mu) =   w(r \tilde g + \mu)$.
Let the dot version of this action be the one shifted by $\rho$
so that $g \cdot \mu =   g(\mu+\rho) -\rho$.

\begin{pro}

(i) For every $g\in G'$, the module $ \Lambda_{g\cdot 0}$ is
${\cal C}_0$-transparent.

(ii) If  $r$ is divisible by $d$,  $r /d >  h^\vee$, and $\zeta$
is a root of unity of order $2Dr$ then every  ${\cal
C}_0$-transparent simple object in ${\cal   C}$  is isomorphic to
$\Lambda_{g \cdot 0}$ with $ g\in G'$.   \label{sub}
\end{pro}

Note that in the expression  $  {g\cdot 0}$ the zero stands for the weight $=
0$.

To prove Proposition \ref {sub} we first recall  the so-called
second symmetry principle for quantum link invariants (see
\cite{Le1}). It describes how  these invariants   change under the
action of $G'$ on the colors of link components. We record two
corollaries of this principle.   First,
\begin{equation} \dim (\Lambda_\nu)^2 =   \dim (\Lambda_{g\cdot \nu})^2,
\quad \forall g\in
G'
\label{dim}
\end{equation}

Second, for the Hopf link $H$,
\begin{equation}
\langle H(\Lambda_{\mu}, \Lambda_{g \cdot \nu}) \,
 \rangle \, \dim(\Lambda_\mu) \, \dim (\Lambda_{g \cdot \nu})=
\ve^{2r (\tilde g|\mu )}\, \langle H(\Lambda_{\mu},
\Lambda_{\nu}) \rangle\, \dim(\Lambda_\mu) \, \dim
(\Lambda_{\nu}) \label{Hopf}
\end{equation}
where $\tilde g\in X'$ is a lift of $g$.

\subsection{Proof of   Proposition \ref {sub}} Part (i).
Set $V=  \Lambda_{g\cdot 0}$. Let $\nu=0$ in both (\ref{dim}) and
(\ref{Hopf}), we see that
$$
\langle H(\Lambda_{\mu}, V)  \rangle = \ve^{2r (\tilde g|\mu )}\,
\dim(\Lambda_\mu) \, \dim (V).$$

Suppose $\mu\in I_0$, then $\mu\in Y$. Hence  $(\tilde g|\mu)\in
\Bbb Z$, since $\tilde g \in X'\subset Y^*$. This implies $\ve^{2r
(\tilde g|\mu )}=1$,  and hence

$$
\langle H(\Lambda_{\mu}, V)  \rangle = \dim(\Lambda_\mu) \, \dim
(V).$$

This means for the Hopf link colored by $V$ and any element in
${\cal C}_0$, we can unlink the Hopf link not altering the value
of the quantum invariant. Thus the operator of Figure \ref{S}a
with $\alpha=0$ is equal to $\Delta_0\,\id$.

Now $\Delta_0= \sum_{\mu \in I_0} \dim(\Lambda_\mu)^2 \neq 0$
since each $\dim(\Lambda_\mu)^2 >0$. By Lemma \ref{trans3}, $V$ is
${\cal C}_0$-transparent.

 Part  (ii).  Consider the colored Hopf link  $H$ depicted in Figure
\ref{last}, where $\omega =   \sum_{\alpha \in G} \omega_\alpha$.
\begin{figure}[htpb]
$$ \psdraw{last}{3.6in} $$
\caption{}\label{last}
\end{figure}

Applying Lemma \ref{trans3} (for strings piercing the unknot with
color $\omega_0$) , we obtain that
$$ \langle H(\omega,\omega_0)\rangle  =   \Delta_0 \sum_{x\in T} (\dim\, x)^
2,$$ where $T$ is the set of isomorphism classes of   objects in
${\cal   M}_{\cal C}$. On the other hand, $\cal C$ is known to be
modular  in the ungraded sense. Hence a string colored by a
simple object  piercing through an unknot with color $\omega$ is
non-zero only when the color of the string is a trivial object.
Hence  $ \langle H(\omega,\omega_0)\rangle =  \langle
U(\omega)\rangle $ which, by Corollary \ref{equal}, is $|G| \,
\Delta_0$. Thus
$$ \sum_{x\in T} (\dim\, x)^2 =   |G|.$$

From Part (i) we know that $\Lambda_{g \cdot 0} \in {{\cal
M}_{\cal C}}$   for every $g\in G'$, and $\dim(\Lambda_{g \cdot
0})^2=1$ by (\ref{dim}). Since $r /d
> h^\vee$, we can have $g \cdot 0=  0$ only when $g=  0$. Hence we
have at least $|G'|=  |G|$ pairwise non-isomorphic   simple
objects in ${\cal M}_{\cal C}$,  each has the  square of the
quantum dimension equal 1. Since  $ (\dim\,   x)^2 $ is a
positive real number for any $x\in T$, we obtain that $T$
consists only of the isomorphism classes of $\Lambda_{g\cdot 0}$
with $g\in G'$.

\begin{rem}
It can be shown  that for every $g\in G'$ one has
$(\Lambda_{g\cdot 0})^*=   \Lambda_{(-g)\cdot 0}$ and  $
\Lambda_{g\cdot 0} \otimes \Lambda_\mu =   \Lambda_{g \cdot \mu}$.
\end{rem}

\subsection{Proof of   Proposition \ref {pre}}
\label{proofofregularity}  We need   to show that ${{\cal M}_{\cal
C}} \subset {\cal C}_0$ and that $\Delta_0^\pm \neq 0$.

Using the explicit description of the weight and root lattices,
one can show  the condition (*) on the number $k$ ensures that

(a)\quad   $r \tilde g \in Y$ for every $\tilde g\in X' = Y^*$,
and

(b) \quad  $d k (\tilde g|\tilde g) \in 2\Bbb Z$ for every
$\tilde g\in X'$.

Now (a) implies that $g \cdot 0 \in Y$ for every $g\in G'$.
Proposition \ref {sub} (ii) shows   that  ${{\cal M}_{\cal C}}
\subset {\cal C}_0$.

Now we prove that $\Delta^\pm_0\neq 0$.  By the second symmetry
principle in \cite{Le1}, the twist $\theta$ acts on  each
$\Lambda_{g\cdot 0}$ as multiplication by  $ \ve^{r d k (\tilde
g|\tilde g)}$  where $\tilde g\in X'$ is a lift of $g\in G'$. Now
(b) implies that  the twist acts as the identity operator on all
${\cal C}_0$-transparent objects.

From (\ref{dim}) we see that for $g \in G'$ one has $
(\dim\Lambda_{g \cdot 0})^2=1$. If $\zeta = \exp(\pi i/Dr)$ (in
that case $\ve = \exp(\pi i/r)$) then $\dim\Lambda_{g \cdot 0}$
is positive, hence $\dim\Lambda_{g \cdot 0}=1$. When $\zeta$ is
an arbitrary root of unity of order $2Dr$, by considering a Galois
action, we also have $\dim\Lambda_{g \cdot 0}=1$.

Note also that $\Delta_0=\sum \dim(x)^2\neq 0$. It follows from
the Bruguieres modularisation criteria  \cite{Bru} (simplified
in  \cite {Blanchet-Beliakova}) that ${\cal C}_0$ is
modularisable. This implies that $\Delta_0^\pm \neq 0$.

\subsection{Proof of Proposition \ref{ref22}} We have $G'= {\Bbb
Z}/2$. Let $g$ be the non-trivial element of $G'$. Explicit
calculation shows that $\Lambda_{g\cdot 0}$ is not in ${\cal
C}_0$. It follows from Proposition \ref{sub} part (ii) that ${\cal
C}_0$ does not have any non-trivial transparent objects. Hence by
a criterion  of Bruguieres (simplified in \cite[Lemma
4.3]{Blanchet-Beliakova}), ${\cal C}_0$ is modular.

\begin{rem} 1. In the $BCD$-case, there is a decomposition of quantum
invariants considered by Blanchet and Beliakova, using
idempotents in the Birman-Wenzl-Murakami category. Their
construction is probably different from ours because they do not
consider spin representations. See also \cite{Sawin}.

2. For a quotient group $\tilde G$ of $G$, one can consider any
$G$-category as a $\tilde G$-category. By choosing $\tilde G$
appropriately, we can make the proof of Proposition \ref{pre}
valid even when the order $r$ of $\ve^2$ does not satisfy the
conditions there. In this way one can get a cohomology
decomposition for $sl_n$ in the case   $1<(r,n)<n$, as in
\cite{Blanchet}. Similar result holds true for the Lie series $D$.

3. Suppose $|G| >1$ (otherwise the theory is reduced to the
ungraded case). If the order $r$ of $\ve^2$ does not satisfy the
conditions of Proposition \ref{ref22} or \ref{ref11}, then it can
be shown that ${\cal C}_0$ contains a non-trivial ${\cal
C}_0$-transparent element of the form $\Lambda_{g\cdot 0}, g\in
G'$, hence $\cal C$ can not be modular (in the $G$-category
sense).
\end{rem}

\end{document}